\documentclass{amsart}

 \usepackage{amsmath}
 \usepackage{amsthm}
\usepackage{amssymb}

\newtheorem{theorem}{Theorem}[section]
\newtheorem{lemma}[theorem]{Lemma}
\newtheorem{proposition}[theorem]{Proposition}
\theoremstyle{definition}

\theoremstyle{remark}
\newtheorem{example}[theorem]{Example}
\newtheorem{remark}[theorem]{Remark}
\numberwithin{equation}{section}
\newtheorem{corollary}{Corollary}[section]

\begin{document}
 \title[Critical points]{On smooth maps with finitely 
 many critical points}
 
\author[D.Andrica]{Dorin Andrica}
\thanks{Partially supported by the European 
 Research Training Network Geometric Analysis.}
\address{Faculty of Mathematics and Computer Science,"Babes-Bolyai"
  University of Cluj, 3400 Cluj-Napoca, \;\;\;\; Romania}
\email{dandrica@math.ubbcluj.ro}

\author[L.Funar]{Louis Funar}
\address{Institut Fourier BP 74, UMR 5582, Universit\'e de  
Grenoble I, 38402 Saint-Martin-d'H\`eres cedex, France}
\email{funar@fourier.ujf-grenoble.fr}

 \begin{abstract} 
We compute the minimum number of critical points 
of a small codimension smooth map between two manifolds. We 
give as well some partial results for the case of higher codimension
when the manifolds are spheres. 
\end{abstract}

\subjclass{57 R 45, 58 K 05, 57 R 60, 14 P 25, 57 R 70}

 
\maketitle
 
 \section{Introduction}
   
If $M, N$ are manifolds, possibly with boundary, consider maps
$f:M\rightarrow N$ with $\partial M=f^{-1}(\partial N)$ such that $f$
has no critical points on $\partial M$. Denote by $\varphi(M,N)$ the
minimal number of critical points of such maps.
The reader may consult the survey \cite{An} for an account 
of various features of this invariant 
(see also \cite{Pi}). Most of the previously known results 
consist of sufficient conditions on $M$ and $N$ 
ensuring that $\varphi(M,N)$ is infinite.

The aim of this note is to find when non-trivial $\varphi(M^m,N^n)$ 
can occur if the dimensions $m$ and $n$ of $M^m$ and respectively 
 $N^n$, satisfy $m\geq n\geq 2$. 
Non-trivial  means here finite and non-zero. 
Our main result is the following: 
\begin{theorem}
Assume that $M^m, N^n$ are compact orientable manifolds and 
$\varphi(M^m, N^n)$ is finite, where 
$0\leq m-n\leq 3$ and $(m,n)\not\in\{(2,2), (4,3), (4,2), (5,3), 
(5,2), (6,3), (8,5)\}$. If $m-n=3$ we also assume 
that the Poincar\'e conjecture in dimension 3 holds true. 

Then $\varphi(M^m, N^n)\in \{0, 1\}$ and $\varphi(M^m, N^n)=1$ 
precisely when the following two conditions are fulfilled: 
\begin{enumerate}
\item  $M^m$ is diffeomorphic to the connected sum 
$\widehat{N}\sharp\Sigma^m$, where $\Sigma^m$ is an exotic sphere 
and $\widehat{N}$ is a $m$-manifold which fibers over 
$N^n$. 
\item  $M^m$ does not fiber over $N^n$.
\end{enumerate}  
\end{theorem}
\begin{proof} 
The statement is a consequence of propositions 3.1, 4.1 and 5.1. 
\end{proof}
\begin{remark}
\begin{enumerate}
\item The second condition is necessary, in general. 
There exist examples 
of connected sums $\widehat{N}\sharp\Sigma^m$
which fiber over $N$, yet
they are not diffeomorphic to $\widehat{N}$. 
In fact, the exotic $7$-spheres constructed 
by  Milnor in \cite{Mi} are pairwise non-diffeomorphic
fibrations over $S^4$ with fiber 
$S^3$.
\item However, if the codimension $m-n$ is zero
we believe that the second condition is redundant
i.e. if $M^m$ is diffeomorphic to 
$\widehat{N}\sharp\Sigma^m$ and $M^m$ is 
not diffeomorphic to $\widehat{N}$ then $M^m$  cannot be a (smooth) 
covering of $N$. This claim holds true when $N^n$ is hyperbolic, 
for all but finitely many coverings $\widehat{N}$. 
In fact, Farrell and Jones (\cite{FJ}) proved that a  
finite covering $\widehat{N}$ of sufficiently large degree
of a hyperbolic manifold $N^n$ has the property that 
$\widehat{N}\sharp \Sigma$ admits a 
Riemannian metric of negative curvature but it does not have 
a hyperbolic
structure. In particular, $\widehat{N}\sharp \Sigma$ is not a covering
of $N$ and hence $\varphi(\widehat{N}\sharp \Sigma, N)=1$.

Conversely if $\varphi(M^n, N^n)=1$  and $N^n$ is
hyperbolic then $M^n$ cannot be hyperbolic. 
Otherwise Mostow rigidity would imply 
that $M^n$ is isometric and hence diffeomorphic to 
$\widehat{N}$. 

\item The theorem holds true for non-compact manifolds $M$ and $N$
if we define $\varphi$ by restricting ourselves 
to those smooth maps which are proper. 

\item In most of the cases excluded in the hypothesis of the theorem 
one can find examples with non-trivial 
$\varphi(M^m, N^n)\geq 2$ (see below). 
\item One expects that for all $(m,n)$ with $m-n\geq 4$ such examples
abound. This is the situation for the local picture. 
The typical example is a  complex projective 
manifold $X$ admitting non-trivial  morphisms into 
${\mathbb C\mathbb P}^1$. 
\item The case $n=1$ was analyzed in \cite{PR}, where the authors 
proved that $\varphi(M,[0,1])=2$, for any non-trivial $h$-cobordism $M$. 
\end{enumerate}
\end{remark}

Most of the present paper is devoted to the proof of the theorem 1.1. 
In the last part  we also compute the values 
of $\varphi(S^m, S^n)$ in a few cases and 
look for a more subtle invariant which would measure how far 
is a manifold from being a covering of another one.

We will consider henceforth that all manifolds are closed  
and connected unless the opposite is stated. 

\vspace{0.3cm}
\noindent {\bf Acknowledgements:} We are grateful to N. A'Campo,
M. Aprodu, E. Ferrand, D.Matei, L.Nicolaescu, C. Pintea, 
P. Popescu-Pampu, A. Sambusetti and J. Seade for useful 
remarks and discussions. 
 
 \section{Elementary computations for surfaces}
Patterson (\cite{Pa}) gave necessary and 
sufficient conditions for the existence of a covering of 
a surface with prescribed degree and ramification orders. 
Specifically his result
 can be stated as follows:
 \begin{proposition}
 Let $X$ be a Riemann surface of genus $g\geq 1$. Let
 $p_1,..., p_k$ be distinct points of $X$ and $m_1,...,m_k$
 strictly positive integers so that
 \[ \sum_{i=1}^k (m_i-1) = 0 \; ({\rm mod} \; 2). \]
 Let $d$ be an integer such that $d\geq max_{i=1,...,k} m_i$.
 Then there exists a Riemann surface $Y$ and a holomorphic
 covering map $f:Y\to X$ of degree $d$ such that there exist
 $k$ points $q_1,...,q_k$ in $Y$ so that $f(q_j)=p_j$,
 $f$ is ramified to order $m_j$ at $q_j$ and is unramified outside the
 set $\{q_1,...,q_k\}$.
 \end{proposition}
 
 \noindent Observe  that a smooth map $f:Y\to X$
 between surfaces has finitely many critical points if and only if
 it is a ramified covering. Furthermore, $\varphi(Y,X)$ is   
 the minimal number of ramification points of a covering $Y\to X$. 
 Estimations can be obtained from the previous result.
 Denote by $\Sigma_g$ the oriented surface of genus $g$.
 Denote by $\left [\left [r\right]\right]$ the smallest 
integer greater than or equal to $r$.
 Our principal result in this section is the following: 
\begin{proposition}
Let $\Sigma$ and $\Sigma'$ be closed oriented surfaces of Euler
characteristics $\chi$ and $\chi'$ respectively. 
\begin{enumerate}
\item If $\chi' > \chi$ then $\varphi(\Sigma',\Sigma)=\infty$. 
\item If $\chi' \leq  0$ then $\varphi(\Sigma', S^2)=3$.
\item If $\chi' \leq -2$ then $\varphi(\Sigma', \Sigma_1)=1$.
\item If $2+2\chi\leq \chi'< \chi \leq -2$ then 
$\varphi(\Sigma',\Sigma)=\infty$. 
\item If $0\leq \mid \chi\mid \leq \frac{\mid\chi'\mid}{2}$, write 
$|\chi'|=a|\chi|+b$ with $0\leq b<|\chi|$; then 
\[ \varphi(\Sigma',\Sigma)=\left[\left[\frac{b}{a-1}\right]\right].\]
\end{enumerate}
In particular, if $G\geq 2(g-1)^2$ then: 
\[ \varphi(\Sigma_G, \Sigma_g)= 
\left \{ \begin{array}{ccc}
 0 & \mbox{\rm if } & \frac{G-1}{g-1}\in {\mathbb Z}_+,\\
 1 & {\mbox {\rm otherwise.}} & 
\end{array}\right.
\]
\end{proposition}
\begin{proof}
The first claim is obvious. 

Further, 
$\varphi(\Sigma', S^2)\leq 3$ because any surface is a covering 
of the 2-sphere branched at three points (from \cite{Al}). 
A deeper result is that the same inequality holds in 
 the holomorphic framework. In fact, Belyi's theorem states that any Riemann surface
 defined over a number field admits a meromorphic function on it with
 only three critical points (see e.g. \cite{SV}).

On the other hand, assume that $f:\Sigma'\to S^2$ is a ramified
covering with at most two critical points. Then, $f$ induces a covering
map $\Sigma'-f^{-1}(E)\to S^2-E$, where $E$ is the set of critical
values and its cardinality $|E|\leq 2$. Therefore one has 
an injective homomorphism $\pi_1(\Sigma'-f^{-1}(E))\to \pi_1(S^2-E)$. 
Now $\pi_1(\Sigma')$ is a quotient of $\pi_1(\Sigma'-f^{-1}(E))$ and 
$\pi_1(S^2-E)$ is either trivial or infinite cyclic, which implies that
$\Sigma'=S^2$.

Next, the unramified coverings of tori are tori; thus any smooth 
map $f:\Sigma_G\to \Sigma_1$ with finitely many critical points 
must be ramified, so that $\varphi(\Sigma_G,\Sigma_1)\geq 1$, if  
$G\geq 2$. 
On the other hand, by Patterson's theorem, there exists a 
covering  $\Sigma'\to \Sigma_1$ of degree 
$d=2G-1$ of the torus, with a single ramification point of 
multiplicity $2G-1$.  From the Hurwitz
formula it follows that $\Sigma'$ has genus $G$, which shows that 
$\varphi(\Sigma_G,\Sigma_1)= 1$. 
 
 \begin{lemma}
 $\varphi(\Sigma', \Sigma)$ is the smallest integer $k$ which satisfies: 
\[  \left[\left[\frac{\chi'-k}{\chi-k}\right]\right] \leq \frac{\chi'+k}{\chi}\]
\end{lemma}
 \begin{proof}
 Suppose that $\Sigma_G$ is a covering of degree $d$ of $\Sigma_g$,
 ramified at $k$ points with the multiplicities $m_i=d-\lambda_i$, where 
 $0\leq \lambda_i\leq d-2$. If one sets $\lambda=\sum_i
 \lambda_i$, then  
$\lambda$ satisfies the obvious inequality:
\[ \lambda \leq k(d-2). \]
Further,  the Hurwitz formula yields the following identity: 
\[ d(k-\chi)=k-\chi'+\lambda.\]
Conversely, if there are solutions $(k,\lambda, d)$ of the two
equations above, with $k, \lambda \geq 0$ and $d\geq 1$, then 
one can find integers $m_i, \lambda_i$ as above and therefore 
one can construct (using 
Patterson's theorem) a ramified covering $\Sigma'\to \Sigma$ of degree
$d$, with $k$ ramification points of multiplicities $m_i$. 
So, $\varphi(\Sigma', \Sigma)$ is the least integer $k\geq 0$ for
which there exists a solution $(k,\lambda, d)\in {\mathbb N}\times {\mathbb
  N}\times {\mathbb N}_+$ of the system:
\[ 0\leq d(k-\chi)+\chi'-k=\lambda \leq k(d-2).\]
That is, for $\chi \leq -2$, $\varphi(\Sigma', \Sigma)$ is the least
$k\in {\mathbb N}$ for which there exists a positive integer $d$ 
satisfying 
\[ \frac{\chi'-k}{\chi-k} \leq d \leq \frac{\chi'+k}{\chi}\]
and this is clearly equivalent to what is claimed in Lemma 2.3. 

%
 \end{proof}
 
 Assume now that   $2+2\chi\leq \chi'< \chi \leq -2$. 
If $f:\Sigma'\to \Sigma$ was a ramified covering then we would have 
$\frac{\chi'+k}{\chi} < 2$, and Lemma 2.3 would imply that
$\chi'=\chi$, which is a contradiction. Therefore 
$\varphi(\Sigma', \Sigma)=\infty$ holds. 
 
Finally, assume that $\frac{\chi'}{2}\leq \chi \leq -2$.
One has to compute the minimal $k$ satisfying 
\[ \left[\left[ \frac{a\chi-b-k}{\chi-k}\right]\right]
\leq \frac{a\chi-b+k}{\chi}, \]
or equivalently, 
\[ \left[\left[ \frac{b+(1-a)k}{\chi-k}\right]\right]
\geq \frac{b-k}{\chi}. \]
The smallest $k$ for which the quantity in the brackets is non-positive 
is $k=\left[\left[\frac{b}{a-1}\right]\right]$, in which case 
\[ \left[\left[ \frac{b+(1-a)k}{\chi-k}\right]\right]\geq 0\geq 
 \frac{b-k}{\chi}. \]
For $k$ smaller than this value one has a strictly positive 
integer on the left-hand side, 
which is therefore at least $1$. But the right-hand side is strictly
smaller than $1$, hence the inequality cannot hold. 
This proves the claim.
 \end{proof}

 \section{Equidimensional case $n\geq 3$}
The situation changes completely in dimensions $n\geq 3$. 
According to (\cite{CH}, II, p.535), H.Hopf was the first to 
notice that  
a smooth map ${\mathbb R}^n\to {\mathbb R}^n$ ($n\geq 3$) which has only an
isolated 
critical point $p$ 
is actually a local homeomorphism at $p$. 
Our result below is an easy application of this fact. 
We outline the proof for the sake of completeness.  
 \begin{proposition}
 Assume that  $M^n$ and $N^n$ are compact manifolds.   
 If $\varphi(M^n,N^n)$ is finite and  $n\geq 3$ 
then $\varphi(M^n,N^n)\in\{0, 1\}$.
Moreover, $\varphi(M^n, N^n)=1$ if and only if $M^n$ is the connected sum 
of a finite covering $\widehat{N^n}$ of $N^n$ with an exotic sphere
and $M^n$ is not a covering of $N^n$.
 \end{proposition}
 \begin{proof}
 There exists a smooth map $f:M^n\to N^n$ which is a local 
 diffeomorphism on the preimage of the complement of a 
 finite subset of points. Notice that $f$ is a proper map. 

 Let $p\in M^n$ be a critical point and $q=f(p)$.   
Let $B \subset N$ be a closed ball intersecting the
set of critical values of
$f$ only at $q$. We suppose moreover that $q$ is an interior
point of $B$. Denote by $U$ the connected component of 
$f^{-1}(B)$ which contains
$p$. As $f$ is proper, its restriction to $f^{-1}(B-\{q\})$ is
also proper. As it is a local diffeomorphism onto $B-\{q\}$, it is a
covering, which implies that $f:U-f^{-1}(q)\rightarrow B-\{q\}$ is also a
covering.
But $f$ has only finitely many critical points in $U$, which
shows that
$f^{-1}(q)$ is discrete outside this finite set, and so
$f^{-1}(q)$ is
countable. This shows that $U-f^{-1}(q)$ is connected. As
$B-\{q\}$ is
simply connected, we see that $f:U-f^{-1}(q)\rightarrow B-\{q\}$
is a
diffeomorphism. This shows that $f^{-1}(q)\cap U=\{p\}$,
otherwise
$H_{n-1}(U-f^{-1}(q))$ would not be free cyclic. So, $f:U-\{p\}
\rightarrow
B-\{q\}$ is a diffeomorphism.
An alternative way is to observe that $f|_{U-\{p\}}$ 
is a proper submersion because $f$ is injective in a neigborhood of $p$ 
(except possibly at $p$). This implies that 
$f:U-\{p\}\to B-\{q\}$ is a covering and hence a diffeomorphism
since $B-\{q\}$ is
simply connected. 

One verifies then easily that the inverse of  
 $f|_{U}:U\to B$ is continuous at $q$ hence it is a homeomorphism. 
 In particular, $U$ is homeomorphic to a ball. Since $\partial U$ is 
 a sphere the results of Smale (e.g. \cite{Sm}) imply that 
 $U$ is diffeomorphic to the ball for $n\neq 4$.   

 We obtained that $f$ is a local homeomorphism hence topologically a 
 covering map. Thus  $M^n$ is homeomorphic to a covering of $N^n$. 
 Let us show now that one can modify $M^n$ by taking the 
 connected sum with an exotic sphere in order to get a smooth covering of $N^n$. 
 
 By gluing a disk to $U$, using an  
 identification 
 $h:\partial U\to \partial B=S^{n-1}$, we obtain a  
 homotopy sphere (possibly exotic) $\Sigma_1=U\cup_{h} B^n$. 
 Set $M_0=M-{\rm int}(U)$, $N_0=N-{\rm int}(B)$.   Given the 
 diffeomorphisms $\alpha:S^{n-1}\to \partial U$ and 
 $\beta:S^{n-1}\to \partial B$ one can form the manifolds 
\[ M(\alpha)=M_0\cup_{\alpha:S^{n-1}\to \partial U} B^n,
N(\beta)=N_0\cup_{\beta:S^{n-1}\to \partial B} B^n.\]
Set $h=f|_{\partial U}:\partial U \to \partial B=S^{n-1}$.  
There is then a map $F: M(\alpha)\to N(h\circ \alpha)$ given by: 
\[ F(x)= \left\{ \begin{array}{ll}
 x & \mbox{ if } x\in D^n, \\
 f(x) & \mbox{ if } x\in M_0.
\end{array}\right. 
\]
The map $F$ has the same critical points as $f|_{M_0}$, hence 
it has precisely one critical point less than $f:M\to N$. 

We choose $\alpha=h^{-1}$ and we remark that  
$M=M(h^{-1})\sharp \Sigma_1$, where the equality sign $"="$ stands 
for diffeomorphism equivalence.  Denote $M_1=M(h^{-1})$. 
We obtained above that  $f: M=M_1\sharp \Sigma_1\to N$ decomposes as follows.   
The restriction of $f$ to $M_0$ extends to $M_1$ without introducing 
extra critical points while the restriction to the homotopy ball 
corresponding to the holed $\Sigma_1$ has precisely one critical
point. 
 
Thus iterating this procedure one finds that there exist possibly exotic
spheres $\Sigma_i$ so that \\ 
$f:M=M_k\sharp \Sigma_1\sharp \Sigma_2...\sharp \Sigma_k\to N$
 decomposes as follows: 
the restriction of $f$ to the $k$-holed $M$ has no critical points, 
and it extends to $M_k$ without introducing any further critical point. 
Each critical point of $f$ corresponds to a (holed) exotic $\Sigma_i$. 
In particular, $M_k$ is a smooth covering of $N$. 

Now  the connected sum $\Sigma=\Sigma_1\sharp \Sigma_2...\sharp \Sigma_k$
is also an exotic sphere.  Let $\Delta=\Sigma-{\rm int}(B^n)$ 
be the homotopy ball obtained by removing an open ball from $\Sigma$.
We claim  that there 
exists a smooth map $\Delta\to B^n$, extending any given diffeomorphism 
of the boundary and  which has exactly one critical point.
Then one builds up a smooth map $M_k\sharp \Sigma\to N$ having 
precisely one critical point, 
by putting together the obvious covering on the 1-holed $M_k$ 
and $\Delta\to B^n$.  This will show that $\varphi(M,N)\leq 1$. 

The claim follows easily from the following two remarks. 
First, the homotopy ball $\Delta$ is diffeomorphic to 
the standard ball by \cite{Sm}, when $n\neq 4$. 
Further, any diffeomorphism $\varphi:S^{n-1}\to S^{n-1}$ 
extends to a smooth homeomorphism with one critical point 
$\Phi:B^n\to B^n$, for example 
\[ \Phi(z)= \exp\left(-\frac{1}{\parallel z\parallel^2}\right)\varphi
\left(\frac{z}{\parallel z\parallel}\right).\]

For $n=4$ we need an extra argument. Each homotopy ball 
$\Delta_i^4=\Sigma_i-{\rm int}(B^4)$ is the preimage $f^{-1}(B)$ 
of a standard ball $B$. Since $f$ is proper we can chose $B$ small 
enough such that  $\Delta_i^4$ is contained in a standard 4-ball. 
Therefore $\Delta^4$ can be engulfed in  $S^4$. 
Moreover, $\Delta^4$ is the closure of one connected component 
of the complement of $\partial \Delta^4=S^3$ in $S^4$. 
The result of Huebsch and Morse from 
\cite{HM} states that any diffeomorphism $S^3\to S^3$ has a 
Schoenflies extension  
to a homeomorphism $\Delta^4\to B^4$ which is a diffeomorphism 
everywhere but at one (critical) point. This proves the claim. 

Remark finally that $\varphi(M^n, N^n)=0$ if and only if 
$M^n$ is a covering of $N^n$. 
Therefore if $M^n$ is diffeomorphic to the 
connected sum $\widehat{N^n}\sharp \Sigma^n$ of a covering $\widehat{N^n}$
with an exotic sphere $\Sigma^n$, and
if it is not diffeomorphic to a covering of ${N^n}$ then 
$\varphi(M^n, N^n)\neq 0$. 
Now drill a small hole  in $\widehat{N^n}$ and glue (differently)
an $n$-disk $B^n$ (respectively a homotopy 4-ball if $n=4$)
in order to get $\widehat{N^n}\sharp \Sigma^n$. The restriction of the
covering $\widehat{N^n}\to N^n$ to the boundary of the hole 
extends (by the previous arguments) to a smooth homeomorphism 
with one critical point over $\Sigma^n$. Thus
$\varphi(M^n, N^n)=1$. 
 \end{proof}
\begin{remark}
\begin{enumerate}
\item We should stress that not all exotic 
structures on a manifold can be obtained from a given 
structure by connected sum with an exotic sphere. For example 
smooth structures on products of spheres (and sphere bundle of spheres) 
are well understood (see \cite{DS1,DS2,Ka,Sch1,Sch2}). 
All smooth structures on $S^p\times S^q$ 
are of the form $(S^p\times \Sigma^q){\#}\Sigma^{p+q}$, where 
$\Sigma^r$ denotes a homotopy $r$-sphere. 
If $p+3\geq q\geq p$ then it is enough to consider only those  
manifolds for which $\Sigma^q=S^q$ (\cite{Ka}) 
but otherwise there are examples where 
the number $n(p,q)$ of non-diffeomorphic manifolds among them 
is larger than the number of homotopy $(p+q)$-spheres. 
For example $n(1,7)=30$, $n(3,10)=4$, $n(1,16)=n(3,14)=24$.

On the other hand, the connected sum with an exotic sphere does not 
necessarily change the diffeomorphism type. 
For example Kreck (\cite{K}) proved that 
for any manifold $M^m$ (of dimension $m\neq 4$)
there exists an integer $r$ such that either 
$M\sharp_r S^m$ 
or $M\sharp ST(S^{\frac{m-1}{2}})\sharp_r S^m$ 
(if $m=1 ({\rm mod }\, 4)$)
has a unique smooth structure, where $ST(S^{k})$ denotes the 
sphere bundle of the tangent bundle of the sphere $S^{k}$. 

However, 
results of Farrell and Jones \cite{FJ} show that any hyperbolic manifold
has finite coverings for which making a connected sum 
with an exotic sphere will change the diffeomorphism type.  
\item  Suppose that $M^n=\widehat{N^n}\sharp \Sigma$  
is not diffeomorphic to  $\widehat{N^n}$. It would be interesting
to know under which hypotheses one can insure that  
$M^n$ is  not a smooth covering of $N^n$.
\end{enumerate}
\end{remark}
 \begin{corollary}
If the dimension $n\in\{3,5,6\}$ then $\varphi(M^n,N^n)$ is either 
$0$ or $\infty$.
\end{corollary}
\begin{proof}
In fact, two 3-manifolds which are homeomorphic are diffeomorphic and  
in dimensions 5 and 6 there are no exotic spheres. 
\end{proof}
\begin{remark}
A careful analysis of open maps between manifolds of the same 
dimensions was carried out in 
\cite{CH}. In particular, one proved that an open map of finite degree 
whose branch locus is a locally tame embedded finite complex 
(for example if the map is simplicial) has both the  branch locus 
and the critical set of codimension 2 (see \cite{CH}, II) 
around each point.   
\end{remark}

\section{Local obstructions for higher codimension}
Our main result in this
section,  less precise than that for codimension $0$, 
is a simple consequence of the 
investigation of local obstructions.  In fact, 
the existence of analytic maps ${\mathbb R}^m\to {\mathbb R}^n$ 
with isolated 
singularities is rather exceptional in the context of smooth real maps
(see \cite{M}). 
\begin{proposition}
If $\varphi(M^m,N^n)$ is finite and either
$m=n+1\neq 4$, $m=n+2\neq 4$, or 
$m=n+3\not\in\{5,6,8\}$ (when one assumes the Poincar\'e conjecture to be
true)  then $M$ is homeomorphic to a fibration of base $N$. 
In particular, if $m=3, n=2$ then 
$\varphi(M^3, N^2)\in \{0, \infty\}$, except possibly for 
$M^3$ a non-trivial homotopy sphere and $N^2=S^2$. 
\end{proposition}
\begin{proof}
One shows first: 
\begin{lemma}\label{alg}
Assume that $\varphi(M^m, N^n)\neq 0$ is finite for two manifolds 
$M^m$ and $N^n$. Then there exists 
a polynomial map $f:({\mathbb R}^m,0)\to ({\mathbb R}^n,0)$ having an isolated
singularity at the origin.
\end{lemma}
\begin{proof}
The hypothesis implies the existence of a smooth map 
$f:({\mathbb R}^m,0)\to ({\mathbb R}^n,0)$ with one isolated singularity 
at the origin. We can assume that the critical point 
is not an isolated point of the fiber over $0$ (see the remark 
\ref{iso} below). 
If the restriction $f|_{S^{m-1}}$ is of maximal rank then the 
construction goes at follows.  
One approximates the restriction  
$f|_{S^{m-1}}$ to the 
unit sphere, up to the first derivative,  
by a polynomial map $\tilde\psi$ (of some degree $d$) and
one extends the later to all of ${\mathbb R}^m$ by 
$\psi(x)=|x|^d \tilde\psi\left(\frac{x}{|x|}\right)$. 
If the approximation is sufficiently close then $\tilde\psi$ will be 
of maximal rank around the unit sphere  hence $\psi$ will have 
an isolated singularity at the origin. 
However, some caution is needed when $f|_{S^{m-1}}$ is not of maximal 
rank. We consider then the restriction $f|_{B^m-B^m_{1-\delta}}$ 
to the annulus bounded by the spheres of radius $1$ and  $1-\delta$ 
respectively. We claim that: 
\begin{lemma}
There exists some $\delta>0$ and
 a polynomial map $\tilde\psi$ (of some degree $d$) 
such that its extension $\psi(x)=|x|^d
\tilde\psi\left(\frac{x}{|x|}\right)$
approximates $f|_{B^m-B^m_{1-\delta}}$ sufficiently close. 
\end{lemma}
\begin{proof}
It suffices to see that $f^{-1}(0)\cap (B^m-B^m_{1-\delta})$ has a
conical structure. Remark that the function 
$r(x)=|x|^2$ has finitely many critical values on 
$f^{-1}(0)\cap (B^m-B^m_{1-\delta})$ since $f$ is smooth and
has no critical points in this range. Thus one can choose 
$\delta$ small enough  so that 
$r$ has no critical points. 
Then the proof of theorem 2.10 
p.18 from \cite{M} applies in this context.
This implies the existence of a good approximation of conical type.  
\end{proof}

In particular, our approximating $\psi(x)$ has no critical points 
in the given annulus. However, if $x_0\neq 0$ was a critical point 
of $\psi(x)$ in the ball, then all points of the line $0x_0$ 
would be critical, since $\psi$ is homogeneous. 
Therefore $\psi$ has an isolated singularity at the origin. 

Notice that one can choose the 
approximation so that $S^{m-1}\cap \psi^{-1}(0)$ is isotopic 
to $S^{m-1}\cap f^{-1}(0)$ and therefore non-empty. 
In particular, the singularities at the origin 
of $f$ and $\psi$  have the same 
topological type.  

Although $\psi$ is not real analytic, each of its components 
are algebraic because they can be represented as  
$\psi_j(x)=P_j(x) + Q_j(x)|x|$, where $P_j(x)$ (resp. $Q_j(x)$) 
are polynomials of even (resp. odd) degree. 
The curve selection lemma (\cite{M}, p.25)  
can be extended without difficulty to sets defined by 
equations like the $\psi_j$ above. Then the proof of theorem 11.2
(and lemma 11.3) from \cite{M}, p.97-99 extends to the case of 
$\psi$. In particular, there exists a Milnor fibration associated to
$\psi$ (the complement of the singular fiber $\psi^{-1}(0)$ in the 
unit sphere $S^{m-1}$ fibers over $S^{n-1}$). Alternatively 
$(S^{m-1}, S^{m-1}\cap \psi^{-1}(0))$ is a Neuwirth-Stallings 
pair  according to \cite{Loo} and $S^{m-1}\cap \psi^{-1}(0)$ is
non-empty. The main theorem from \cite{Loo} 
provides then a polynomial map with an isolated 
singularity at the origin as required.  
\end{proof}

Milnor (see \cite{M}) called such an isolated singularity 
{\em trivial} when its local Milnor fiber is diffeomorphic to a disk. 
Then it was shown in (\cite{CL}, p.151) that $f$ is trivial if and 
only if $f$ is locally topologically equivalent to the projection map 
${\mathbb R}^m \to {\mathbb R}^n$, whenever the dimension of the fiber is 
$m-n\neq 4, 5$.  We recall that the existence of polynomials 
with isolated singularities was (almost) settled in \cite{CL,M}: 
\begin{proposition}\label{cl}
For $0\leq m-n \leq 2$ non-trivial polynomial singularities 
exist precisely for $(2,2)$, $(4,3)$ and $(4, 2)$. 

For $m-n \geq 4$ non-trivial examples occur for all $(m,n)$. 

For $m-n=3$ non-trivial examples occur for $(5,2)$ and $(8,5)$. 
Moreover, if the 3-dimensional Poincar\'e conjecture is false 
then there are non-trivial examples for all $(m,n)$. Otherwise 
all examples are trivial except $(5, 2)$, $(8, 5)$ and possibly 
$(6, 3)$.    
\end{proposition}
We consider now a smooth map $f:M^m\to N^n$ where $m, n$ 
are as in the hypothesis. For each critical point $p$ there 
are open balls $2B^m(p)$ and $2B^n(f(p))$ for which the restriction 
$f|_{2B^m(p)}:2B^m(p)\to 2B^n(f(p))$ has an isolated singularity at $p$.
One identifies  $2B^m(p)$ 
with the ball of radius 2 in ${\mathbb
  R}^m$, and let $B^m(p)$ 
be the preimage of the concentric unit ball.  
In the proof  of lemma \ref{alg} we  
approximated $f|_{\partial B^m(p)}$ by a polynomial map $g$
with isolated singularities, both maps having isotopic links and 
being close to each other. 
Assume for simplicity that $f$ (hence $g$) is 
of maximal rank around  this (unit) sphere. The general case follows 
along the same lines. 
Then there exists an isotopy  $f_t$ $(t\in [0,1])$ between 
$f|_{\partial B^m(p)}$ and $g|_{\partial B^m(p)}$, which is close 
to identity. In particular, all $f_t$ are of maximal rank around
the unit sphere. Let $\rho:[0,4]\to [0,1]$ be a smooth  
decreasing function with 
$\rho(x)=0$ if  $x\geq 1$ and 
$\rho(x)=1$ if $x\leq \frac{1}{2}$. 
Let $F:2B^m(p)\to 2B^n(f(p))$ be the map defined by: 
\[ F(x)= \left\{ \begin{array}{ll}
 f_{\rho(|x|^2)}\left(\frac{x}{|x|}\right) & \mbox{ if } |x|\geq \frac{1}{2}, \\
 g(x) & \mbox{ if } |x|\leq \frac{1}{2}.
\end{array}\right. 
\]
If one replaces $f_{2B^m(p)}$ by $F$ then 
one obtains a smooth function with an 
isolated singularity at $p$, which must be a topological submersion at
$p$ (by the previous proposition \ref{cl}). 
An induction on the number of critical points yields a map 
$F:M^m\to N^n$ which is a topological fibration. 
\end{proof}
\begin{remark}
Notice that there exist real smooth maps $f$ which 
don't have a Milnor fibration at an isolated singularity. 
For such $f$ it is not clear when one should call the singularity 
trivial. In particular, in this situation we don't know 
whether $f$ itself must be a topological
submersion. Therefore it is necessary to replace $f$  by another 
map (locally algebraic), 
in order to be able to apply proposition \ref{cl}. 
\end{remark}

\begin{remark}
Therefore, within the range $0\leq m-n \leq 3$, with the exception 
of $(2,2)$, $(4,3)$, $(4, 2)$, $(5, 2)$, $(8, 5)$ and 
$(6, 3)$ the non-triviality of $\varphi$ is related to the exotic
structures on fibrations.
\end{remark}

One expects that in the case when 
non-trivial singularities can occur  such   
examples abound.

\begin{example} In the remaining cases we have: 

$(m,n)\in \{(4,3), (8,5)\}$.  
We will prove below that $\varphi(S^4, S^3)=\varphi(S^8, S^5)=2$. 

$(m,n)=(4,2)$. Non-trivial examples come from  
Lefschetz fibrations  $X$ over a Riemann surface $F$. 
For instance  $X$ is an elliptic K3 surface and $F$ is 
$\mathbb C\mathbb P^1$. 

$(m,n)=(2k, 2)$. More generally, 
one can consider complex projective $k$-manifolds 
admitting morphisms onto an algebraic curve. 
\end{example}

Further, one notices that these local obstructions are far from being 
complete. In fact, the maps ${\mathbb R}^m\to {\mathbb R}^n$
arising as restrictions of smooth maps between 
compact manifolds are quite particular. 
For instance if one takes $M=S^m$ then one can obtain 
by restriction a map ${\mathbb R}^m\to {\mathbb R}^n$ which is 
proper and has only finitely many isolated singularities. However, 
adding extra conditions can further restrict the range of 
dimensions: 
\begin{proposition}\label{prope}
There are no proper smooth functions 
$f:({\mathbb R}^m,0)\to ({\mathbb R}^n,0)$ with one isolated singularity 
at the origin if $m \leq 2n-3$. 
\end{proposition}
\begin{proof} 
There is a direct proof similar to that of proposition \ref{spcon}. 
Instead, let us show that the hypothesis implies that 
$\varphi(S^m, S^n)\leq 2$ and so proposition \ref{spcon} yields the 
result. 

Let $j_k:S^k\to {\mathbb R}^k$ denote  the stereographic 
projection from the north pole ${\infty}$. 
There exists an increasing unbounded real function $\rho$ such that 
$|f(x)|\geq \lambda(|x|)$ for all $x\in {\mathbb R}^m$, because $f$ is
proper. 

We claim that there exists a real function $\rho$ such that 
$\rho(|x|)f(x)$ extends to a smooth function $F:S^m\to S^n$. 
Specifically, we want that the function $F_{\rho}:S^m\to S^n$ defined by: 
\[ F_{\rho}(x)= \left\{ \begin{array}{ll}
 j_n^{-1}(\rho(|j_m(x)|)f(j_m(x)) & \mbox{ if } x\in S^m-\{\infty\}, \\
 \infty & \mbox{otherwise }.
\end{array}\right. 
\]
be smooth at $\infty$. This is easy to achieve by taking 
$\rho(x) > \exp(|x|) \lambda^{-1}(|x|)$ for large $|x|$. 
Now the critical points of $F_{\rho}$ consist of the two poles, 
and the claim is proved.   
\end{proof}
 \begin{remark}
Notice that proper maps like above for $m=2n-2$ exist only for 
$n\in \{2,3,5,9\}$ (see below). 
\end{remark}
  
\begin{remark}\label{iso}
A special case is when the critical point $p$ is an isolated point 
in the fiber $f^{-1}(f(p))$. This situation was settled in   
\cite{CL,Ti} where it was shown that the dimensions 
$(m, k)$ should be $(2,2), (4,3), (8,5)$ or $(16, 9)$, and 
the map is locally the cone over the respective Hopf fibration. 
\end{remark}

\section{The global structure for topological submersions}
Roughly speaking the results of the previous section 
say that 
maps of low codimension with only 
finitely many critical points should be topological submersions. 
\begin{proposition}
Assume that there exists a  topological submersion $f:M^m\to N^n$ 
with finitely many critical points, and $m>n\geq 2$. 
 Then $\varphi(M,N)\in \{0, 1\}$ and 
$\varphi(M,N)=1$ precisely when $M$ is diffeomorphic 
to the connected sum of a fibration $\widehat{N}$ (over $N$) 
with an exotic sphere, and $M$ is not a fibration over $N$. 
\end{proposition}
\begin{proof}
The first step is to split off one critical point by localizing it 
within an exotic sphere. Let $M_0$ be the manifold obtained after 
excising an embedded ball from $M$. 

\begin{lemma}
There exists an exotic sphere $\Sigma_1$ and a map 
$f_1:M_1\to N$ such that:
\begin{enumerate}
\item $M_0$ is a submanifold of both $M$ and $M_1$.
The complements $M_1-M_0$ and $M-M_0$ are balls and  
$M=M_1\sharp \Sigma_1$. 
\item $f_1$ agrees with $f$ on $M_0\subset M_1$ and has no other 
critical points on the ball $M_1-M_0$.  
\item $f$ has precisely one critical point in $M-M_0$.
\end{enumerate}
\end{lemma}
\begin{proof}
Let $p$ be a critical point of $f$, $q=f(p)$, $\delta$ be a small 
disk around $q$. We replace $f$ by a map which is locally polynomial 
around the critical point $p$, as in the previous section. 
We show first that: 
\begin{lemma}\label{cyl}
There exists a neighborhood $Z_p$ of $p$ such that the following 
conditions are fulfilled: 
\begin{enumerate}
\item $Z_p$ is diffeomorphic to $D^n\times D^{m-n}$ (for $m\neq 4$). 
\item $\partial Z_p=\partial^h Z \cup \partial^v Z_p$, where the
  restrictions $f:\partial^vZ_p\to D^n$ and  
$f:\partial^hZ\to \partial D^n$ are trivial fibrations, and 
$\partial^h Z \cap \partial^v Z_p=S^{n-1}\times S^{m-n-1}$. 
\end{enumerate}
\end{lemma}
\begin{proof}
Let $B^m(p)$ be a sufficiently small ball around $p$ and $\delta$ 
be such that $\delta\subset f(B^m(p))$. 
We claim that $Z_p=B^m(p)\cap f^{-1}(\delta)$ has the required 
properties. 

One chooses a small ball containing $p$, $B_0^m(p)\subset Z_p$. 
Then one uses the argument from (\cite{M}, p.97-98) and derive that  
$q$ is a regular value of the map  
\[ f: Z_p-{\rm int}(B_0^m(p))  \to \delta. \]
Therefore the latter is a fibration, hence a trivial fibration. 
In particular, the manifold with corners
$\partial Z_p$ has a collar whose outer boundary 
is a smooth sphere. 
Further, the manifold with boundary $Z_p$ is homeomorphic 
to $D^n\times D^{m-n}$ and the boundary $\partial Z_p$ is collared 
as above. The outer sphere 
bounds a smooth disk (by Smale) and so $Z_p$ is diffeomorphic
to  $D^n\times D^{m-n}$. 
\end{proof}

Now the proof goes on as in codimension 0. We excise $Z_p$ and glue it
back by another diffeomorphism in order that the restriction of $f$ 
extends over the new ball, without introducing  any new critical points. 
The gluing diffeomorphism respects the corner manifold structure. 
\end{proof}

\noindent An inductive argument shows that if $\varphi(M,N)$ is finite
then the connected sum $M\sharp \Sigma$ with an exotic sphere 
is diffeomorphic to a fibration over $N$. 

We want to prove now that one can find another map $M\to N$ having 
precisely one critical point. We have first to put all critical 
points together inside a standard neighborhood: 
\begin{lemma}
If $m>n\geq 2$ then the critical points of $f$ are 
contained in some cylinder  $Z^m \subset M$ 
which is diffeomorphic to $D^n\times D^{m-n}$ (respectively 
homeomorphic when $m=4$, by a homeomorphism  which is a 
diffeomorphism on the boundary) such that the fibers of $f$ 
are either transversal to the boundary 
(actually to the part $D^{n}\times \partial D^{m-n}$) or contained in
$\partial D^n\times D^{m-n}$. 
\end{lemma}
\begin{proof}
Pick-up a regular point $x_0$ in $M$. Let $U$ be the set
of regular points which can be joined to $x_0$ by an arc $\gamma$ 
everywhere transversal to the fibers of $f$ (which will be called 
transversal in the sequel). 

We show first that $U$ is open. In a small neighborhood $V$ 
of $x\in U$ the fibers can be linearized (by means of a
diffeomorphism) and identified to parallel $(m-n)$-planes. 
Let $y\in V$. If $x$ and  $y$ are not in the same
fiber then the line joining them is a transversal arc. 
Otherwise use a helicoidal arc spinning around the line, which 
can be constructed  since the fibers have codimension at least 2.

At the same time  $U$ is closed in the complement of the critical set. 
In fact, the previous arguments show that 
two regular points which are sufficiently closed to each other 
can be joined by a transversal arc  
with prescribed initial velocity (provided this 
tangent vector is also transversal to the fiber).  
Thus, if $y_i$ converge to a regular point $y$ and $y_i\in U$ then 
$y$ can be joined to $x_0$ by joining first $x_0$ to $y_i$ and further 
$y_i$ to $y$ (for large enough $i$) with some prescribed initial
velocity, in order to insure the smoothness of the arc.  
This proves that $U$ is the set of all regular points. 

Further, we consider the cylinders $Z_{p_i}$ given by lemma \ref{cyl}. 
Let $f_i\subset \partial Z_{p_i}$ be some fibers in the boundary. 
The points $q_i\in f_i$ can be joined by everywhere transversal 
arcs. Since this is an open condition one can find disjoint tubes 
$T_{i, i+1}$ joining neighborhoods of the fibers $f_i$ in $\partial Z_{p_i}$ 
and  $f_{i+1}$ in $\partial Z_{p_{i+1}}$ and one builds up 
this way a cylinder $Z$ containing all critical points.  
\end{proof}

\begin{lemma}\label{one}
There exists a smooth map $g:Z\to D^n$ having one critical point 
such that $g|_{\partial Z}=f|_{\partial Z}$. 
\end{lemma}
\begin{proof}
We know that each restriction  $f:Z_i\to D^n$ is topologically a 
trivial fibration, whose restriction to a collar of the boundary 
sphere $\partial D^n$ is trivial, as a smooth fibration.
We claim that $f:Z\to D^n$ enjoys the same property. This follows from
the fact that the restrictions $f:T_{i, i+1}\to D^n$  
are also trivial fiber bundles. 

The restriction $f|_{\partial^v Z}:\partial^v Z\to D^n$ is a trivial fibration 
over the ball. Thus there exists a diffeomorphism 
$\partial^v f:\partial^v Z\to D^n\times \partial D^{m-n}$, which 
commutes with the trivial projection $\pi:D^n \times D^{m-n}\to D^n$, 
namely $\pi \circ \partial^v f = f|_{\partial^v Z}$.

Since $f|_{\partial^h Z}:\partial^v Z\to \partial D^n$
is a trivial fibration there exists a diffeomorphism 
$\partial^h f:\partial^v Z \to \partial D^n\times D^{m-n}$ commuting with
$\pi$. Moreover, these two diffeomorphisms can be chosen to agree 
on their common intersection, 
$\partial^h f|_{\partial^h Z \cap \partial^v Z}=
\partial^v f|_{\partial^h Z \cap \partial^v Z}$. 

We obtain therefore a diffeomorphism $\partial f$ of manifolds with corners
$\partial f: \partial Z\to \partial D^m$, defined by: 
\[ \partial f(x)= \left\{ \begin{array}{ll}
\partial^h f(x)  & \mbox{ if } x\in \partial^h Z, \\
\partial^v f(x)  &  \mbox{ if } x\in \partial^v Z .
\end{array}\right. 
\]

Assume now that there exists a smooth homeomorphism  
$\Phi:Z\to D^n\times D^{m-n}$ having precisely one critical point,
which extends $\partial f$, i.e. such that the diagram 

\[\begin{array}{ccc} 
\partial Z & \stackrel{\partial f}{\to} & \partial(D^n \times D^{m-n})
\\
\downarrow & & \downarrow \\
Z & \stackrel{\Phi}{\to} & D^n \times D^{m-n} 
\end{array}
\]
commutes. We set therefore 
$g(x)= \pi(\Phi(x))$. It is immediate that $g$ has at most 
one critical point and $g$ is an extension of  $f_{\partial Z}$. 
Our claim is then a consequence of the following: 
\begin{lemma}
Any diffeomorphism of the sphere $S^m$, with the structure of a
manifold with corners $\partial(D^n\times D^{m-n})$, 
extends to a smooth homeomorphism of $D^n\times D^{m-n}$ with at most one 
critical point. 
\end{lemma}
\begin{proof}
Instead of searching  for a direct proof 
remark that the trivializations 
leading to $\partial f$ extend over a collar of $\partial(D^n\times
D^{m-n})$. This collar is still a manifold with corners, but it
contains a smoothly embedded sphere. We use then the standard 
result (see \cite{HM}) to extend further the diffeomorphism from 
the smooth sphere to the ball.    
\end{proof}
\noindent Now lemma \ref{one} follows. 
\end{proof}

\end{proof}

\section{Maps between spheres}
For spheres the situation is somewhat simpler than in general  
because we can use the global obstructions of topological nature. 
Our main result settles the case when the codimension is smaller than 
the dimension of the base. Specifically, we can state: 

\begin{proposition}\label{spcon} 
\begin{enumerate}
\item The values of $m>n>1$  for which $\varphi(S^m,S^n)=0$ 
are exactly those arising in the Hopf fibrations i.e. $n\in\{2,4,8\}$ and 
$m=2n-1$.

\item One has 
$ \varphi(S^4, S^3)=\varphi(S^8, S^5)=\varphi(S^{16}, S^9)=
2$.

\item If $m\leq 2n-3$ then $\varphi(S^m,S^n)=\infty$. 

\item If $\varphi(S^{2n-2},S^n)$ is finite  then $n\in\{2,3,5,9\}$. 
\end{enumerate}
\end{proposition}
\begin{proof}
Notice first that the existence of the Hopf fibrations 
$S^3 \to S^2,  S^7 \to S^4, S^{15}\to S^8$, 
shows that 
$\varphi(S^3, S^2)=\varphi(S^7, S^4)=\varphi(S^{15}, S^8)=0.$
The converse is 
already  known (see Lemma 2.7 in \cite{Ti} or Lemma 1 in 
\cite{CL}). We will give a slightly different proof below, on 
elementary grounds. 
 
Using the Serre exact sequence for a $(n-1)$-connected basis
one finds that the homology of the fiber $F$ 
(of the fibration $f:S^m\to S^n$) agrees with that of $S^{n-1}$ 
up to dimension $n-1$, and a subsequent application of the same 
sequence shows that $F$ is a homology $(n-1)$-sphere. 
Therefore $m=2n-1$. 
In particular, we infer that  the transgression map 
$\tau^*:H^{n-1}(F) \to H^n(S^n)$ is an isomorphism.

Let $i_F:F\to S^{2n-1}$ and 
$j_F:S^{2n-1}\to (S^{2n-1},F)$ denote the inclusion maps.
Denote by $C^*(X)$ the cochain complex of the space $X$.  

\begin{lemma}
The composition of maps:
\[ C^{n-1}(S^{2n-1}) \stackrel{i_F^*}{\to}
  C^{n-1}(F)\stackrel{\tau^*}{\to}
C^{n}(S^{2n})\stackrel{f^*}{\to}
C^{n}(S^{2n-1})\]
is the boundary operator $d:C^{n-1}(S^{2n-1})\to 
C^{n}(S^{2n-1})$. 
\end{lemma}
\begin{proof}
The transgression map $\tau^*$ can be identified (see for example
\cite{Wh}, p.648-651) with the composition 
\[ H^{n-1}(F) \stackrel{\partial^*}{\to} H^n(S^{2n-1},F)
\stackrel{f^{*-1}}{\to} H^n(S^n),\]
where $\partial^*$ is the boundary homomorphism in the 
long exact sequence of the pair $(S^{2n-1}, F)$. 
One sees then that  
\[ C^n(S^{2n-1},
F)\stackrel{f^{*-1}}{\to}C^n(S^n)\stackrel{f^{*}}{\to}
 C^n(S^{2n-1})\]
agrees with  $j_F^*$. 
Further,  the composition of maps from the statement of the lemma is 
equivalent to 
\[  C^{n-1}(S^{2n-1}) \stackrel{i_F^*}{\to}
  C^{n-1}(F)\stackrel{\partial^*}{\to}
C^{n}(S^{2n-1},F)\stackrel{j_F^*}{\to}
C^{n}(S^{2n-1}), \]
which acts as the boundary operator $d$, as claimed.
\end{proof}

Let $u$ be an $(n-1)$-form on 
$S^{2n-1}$ such that $i_F^*u$ is a generator for 
$H^{n-1}(F,{\mathbb Z})\subset H^{n-1}(F, {\mathbb R})$. Then 
\[ <i_F^*u, [F]> = \int_F u =1. \]

Since $\tau^*$ is an isomorphism it follows that 
$\tau^*i_F^*u=v$, where $v$ is the generator of 
$H^{n}(S^n, {\mathbb Z})$. Thus 
$v$ is the volume form on $S^n$, normalized so that 
$\int_{S^n}v=1$. 

Let us recall the definition of the Hopf invariant $H(f)$. 
Consider any $(n-1)$-form $w$ on $S^{2n-1}$ satisfying $f^*v=dw$.  Then 
\[ H(f) = \int_{S^{2n-1}}w\wedge dw =\int_{x\in S^n} \left(\int_{f^{-1}(x)} w \right)  v.\]

According to the lemma one has $f^*v=du$. However, it is clear that
the function $x\to  \int_{f^{-1}(x)}w$ is constant (more generally 
it is locally constant on the set of regular values for an arbitrary 
$f$), and this constant in our case is 
$\int_Fu=1$. Therefore $H(f)=1$ and the  Adams theorem 
(see \cite{Ad}) implies the claim.

\begin{remark}\label{serre}
The result above holds true if one relaxes the assumptions by   
asking $f$ to be only a Serre fibration. One replaces in the proof 
the integral in the definition of the Hopf invariant 
by the intersection of chains
(see e.g. \cite{Wh}, p.509-510).
\end{remark}

\noindent {\em Proof of 2.} We will show now that by suspending the Hopf 
fibrations we obtain  examples of 
pairs with non-trivial $\varphi$. 
In fact, choose a Hopf map $f:S^{2n-1}\to S^n$, and extend it 
to $B^{2n}\to B^{n+1}$ by taking the cone  and smoothing it 
at the origin. Then glue together two copies of $B^{2n}$ 
along the boundary. 
One gets a smooth map having two critical points. 
The previous result implies that: 
\[ 1\leq \varphi(S^4, S^3), \varphi(S^8, S^5), \varphi(S^{16}, S^9)\leq 2.\]

 Let us introduce some more notations: set $p_1,...,p_r$ for the 
critical points of the map $f:S^m\to S^n$ under consideration, if
there are only finitely many. Let $F_{e_i}=f^{-1}(f(p_i))$ denote the 
singular fibers, $F_e=\cup_{i=1}^r F_{e_i}$ stand for their union, 
and $F$ for the generic fiber which is a closed oriented $(m-n)$-manifold.

\begin{lemma}
Each component of $F_e$ is either a smooth $(m-n)$-manifold 
around each point which is not in the critical set $\{p_1,...,p_r\}$, 
or else an isolated $p_i$.  
\end{lemma}
\begin{proof}
In fact, $f$ is a submersion at all points but $p_i$. 
\end{proof}

\begin{lemma}
If $m < 2n-1$ then $\varphi(S^m, S^n)\geq 2$. 
\end{lemma}
\begin{proof}
Assume  that there is a map $f:S^m\to S^n$ with 
precisely one critical point $p$. 
Then  $f:S^m-F_e\to B^n$  is a fibration, so that
$S^m-F_e=B^n\times F$.

One rules out the case when the exceptional fiber is one point
by observing that $H_{m-n}(F)$ is not trivial. 
Using an $(n-1)$-cycle linking once a component of $F_{e_i}$ 
one shows that  $H_{n-1}(S^m-F_e)$ is non-trivial. Since $n-1> m-n$ the 
equality above is impossible, and the claim is proved.  
\end{proof}
\noindent Now the equalities from the statement follow.

\noindent 
\begin{remark}
This might be used to construct other examples with finite 
$\varphi$ in the respective dimensions. For instance one finds that 
$\varphi(\Sigma^8, S^5)=\varphi(\Sigma^{16}, \Sigma^9)=2$, where $\Sigma^n$ 
denotes an exotic $n$-sphere.
\end{remark}

\noindent {\em Proof of 3.}
Assume that there is a smooth map $f:S^m\to S^n$ with 
$r$ critical points. We  suppose, for simplicity, that 
the critical values $q_i$ are distinct.
One uses the Serre exact sequence for the fibration $S^m-F_e\to
S^n-\{q_1,...,q_r\}$ and one derives that: 
\[ H_i(F)=H_i(S^m-F_e), \mbox{ if } i\leq n-3.\]
Further,  $H^{m-i}(F_e)=0$ for $i\leq n-1$, because 
$F_e$ has dimension at most $(m-n)$. Then Alexander's duality, 
$\widetilde{H}_{i-1}(S^m-F_e)=\widetilde{H}^{m-i}(F_e)$, and the previous 
equality imply that  $H_i(F)=0$, for all  $i\leq n-3$. 
This is impossible because  the fiber $F$ is a
compact $(m-n)$-manifold and  $m-n\leq n-3$.

\noindent {\em Proof of 4.}
As above, the Serre exact sequence shows that  $F$ is 
an $(n-2)$-homology sphere. Further, the generalized Gysin sequence 
yields: 
\[ \widetilde{H}^{2n-3-j}(F_e)=\widetilde{H}_j(S^{2n-2}-F_e)=\left\{
\begin{array}{ll}
{\mathbb Z}^{r} & \mbox{ if } j=2n-3, \\
{\mathbb Z}^{r-2} & \mbox{ if } j=n-1, \\
0 & \mbox{ otherwise.}
\end{array}
\right.  
\]
Notice that $H_{n-2}(F_e)$ (or equivalently 
$H_{n-1}(S^{2n-2}-F_e)$) cannot be of rank $r-2$ unless some 
(more precisely two such) exceptional fiber  in $F_e$ consists 
of one point. In fact, if we have $q$ connected components of $F_e$  
of dimension $(n-2)$ then the rank of $H_{n-2}(F_e)$ would be at least
$q$. 

Furthermore, such a critical point $p$ is isolated in $f^{-1}(f(p))$. 
Then proposition 3.1 from \cite{Ti} yields the claim. 
\end{proof}

\begin{remark}
Notice that $F^{m-n}$ is $(n-3)$-connected. In particular, if 
$m\leq 3n-6$ ($n\geq 5$) then $F$ is homeomorphic to $S^{m-n}$. 
In fact, one can obtain  $S^m$ from the complement 
$S^m-{\rm int}(N(F_e))$ of a neighborhood of the exceptional 
fibers  by adding cells of
dimension $\geq n$, one $(n+i)$-cell for each $i$-cell of $F_e$. 
Therefore $\pi_j(S^m-{\rm int}(N(F_e)))=0=\pi_j(S^m)=0$, for $j\leq n-2$. 
The base space of the fibration $f|_{S^m-{\rm int}(N(F_e))}$ is 
$S^n$ with small open neighborhoods of the critical values deleted, 
thus it is homotopy equivalent to a bouquet of $S^{n-1}$ (at least one
critical value). The long exact sequence in homotopy shows then that
the fiber is $(n-3)$-connected. 
\end{remark}

\begin{remark}
\begin{enumerate}
\item If there exists a non-trivial 
proper smooth 
$F:({\mathbb R}^m,0)\to ({\mathbb R}^n,0)$ having only one isolated
singularity at the origin   
then $\varphi(S^m, S^n)\leq 2$ 
(see the proof of proposition \ref{prope}), 
and this condition  seems to be 
quite restrictive, in view of proposition \ref{spcon}. 
\item The explicit computation of $\varphi(S^m,S^n)$  for general $m, n$ 
seems to be difficult. Further steps towards the answer 
would be to prove that 
$\varphi(S^{2n-1},S^n)$ is finite 
only if $n\in\{1,2,4,8\}$, and that  $\varphi(S^m,S^n)=\infty$
if $n\geq 5$ and $2n-1 < m \leq 3n-6$. 
\end{enumerate}
\end{remark}

\section{Remarks concerning a substitute for $\varphi$ in dimension 3}
 One saw that $\varphi(M^n,N^n)$ is less interesting if $n\geq
 3$. One would like to have an invariant of the pair
 $(M^n,N^n)$ measuring how far 
 is $M^n$ from being an unramified covering of $N^n$. First, one has
 to know whether there is a branched covering $M^n\to N^n$ and next
 if the branch locus could be  empty. 
\begin{remark}
 A classical 
 theorem of Alexander (\cite{Al}) states that any $n$-manifold is a branched 
 covering of the sphere $S^n$. Moreover, one can assume that 
the ramification locus is the $(n-2)$-skeleton of the 
standard $n$-simplex. 
\end{remark}
\begin{remark}
There exists 
an obvious obstruction to the existence of a ramified covering 
$M^n\to N^n$, namely the existence of a map of non-zero degree $M^n\to N^n$.
In particular,  a necessary condition is
$\parallel M\parallel\;\geq\; \parallel N \parallel$, 
where $\parallel M\parallel$ denotes the simplicial volume of $M$
 (see \cite{Gro,MT}). 
However, this condition is far from being sufficient. 
Take $M$ with finite fundamental group and $N$ with infinite amenable
fundamental group (for instance of polynomial growth); then 
$\parallel M\parallel=\parallel N \parallel=0$, while it is elementary
that there does not exist a non-zero degree map $M\to N$. 
\end{remark}

 A possible candidate for replacing $\varphi$ in dimension 3
 is the ratio of simplicial volumes mod ${\mathbb
    Z}$, namely
 \[ v(M,N) = \frac{\parallel M\parallel}{\parallel N \parallel}
 \; ({\rm mod}  \; {\mathbb Z})\in [0,1), \]
 which is defined when $N$ has
 nonzero simplicial volume.
 Notice that for (closed manifolds $M$) the simplicial
 volume  $\parallel M\parallel$ depends only on the fundamental group
 $\pi_1(M)$ of $M$. In particular, it vanishes for simply connected
 manifolds, making it less useful in dimensions at least 4.

 \begin{remark}
 If $M^n$ covers $N^n$ then 
 $v(M,N)=0$ (see \cite{Gro}). The converse holds true for surfaces
 of genus at least 2, from Hurwitz formula.
 \end{remark}

 The norm ratio  has been extensively studied
 for hyperbolic manifolds  in dimension 3, where it coincides with
 the volume ratio, in connection with commensurability problems
(see e.g. \cite{T}).  In particular, the values $v(M^3,N^3)$ accumulate on 1
since the set of volumes of closed hyperbolic 3-manifolds has an accumulation
point. The simplicial volume is zero for a Haken 3-manifold iff the manifold 
is a graph manifold (from \cite{So}), and conjecturally the 
simplicial volume is the sum of (the hyperbolic) volumes of the 
hyperbolic components of the manifold.

 However, it seems that this invariant is not appropriate in dimensions
 higher than 3 (even if one restricts to aspherical manifolds). 
 Here are two arguments in  favour of this claim: 
 
 \begin{proposition}
Let us suppose that $M^n$ is a ramified covering of $N^n$ over
 the complex $K^{n-2}$. Assume that both the branch locus $K^{n-2}$ and its
 preimage in $M^n$ can be
 engulfed in a simply connected codimension one submanifold.
 Then $v(M^n,N^n)=0$.

Assume that there is a map $f:M^n\to N^n$ such that the kernel 
$\ker(f_*:\pi_1(M)\to \pi_1(N))$ is an amenable group. 
Then $v(M^n,N^n)=0$. 
\end{proposition}
 \begin{proof}
 One uses the fact that for any simply connected codimension one submanifold
 $A^{n-1}\subset M^n$ one has
 $\parallel M\parallel =\parallel M-A \parallel$ (see \cite{Gro},
 p.10 and 3.5). The second part follows from (\cite{S}, Remark 3.5) 
which states that
$\parallel M\parallel={\rm deg}(f) \parallel N \parallel$, 
where ${\rm deg}(f)$ states for the degree of $f$.
\end{proof}
\begin{remark}
\begin{enumerate}
\item For  all $n\geq 4$ Sambusetti (\cite{S}) constructed examples of 
manifolds $M^n$ and $N^n$ 
satisfying the second condition (and hence $v$ vanishes) 
but which are not fibrations. 
\item It seems that there are no such examples in dimension 3.
At least for Haken hyperbolic $N^3$, any $M^3$ dominating $N^3$ with
amenable kernel must be a covering, according to (\cite{S}, Remark 3.5)
and to the rigidity result of Soma and Thurston (see \cite{So}).   
\end{enumerate}
\end{remark}

\begin{remark} 
One could replace the simplicial volume by any other volume, as 
defined by Reznikov \cite{R}. For instance the $\widetilde{SL(2,R)}$-volume 
is defined for Seifert fibered 3-manifolds and it behaves
multiplicatively under finite coverings (compare to \cite{WaWu}). 
In particular, one can define an appropiate $v(M,N)$ for graph
manifolds using this volume. Other topological invariants which behave
multiplicatively under finite coverings are the $l^2$-Betti numbers. 
\end{remark}
\begin{remark}
 If there is a branched covering $f:M^n\to N^n$ then 
 the branch locus is of codimension 2. This yields a heuristical 
explanation for the almost triviality of $\varphi(M^n,N^n)$ 
in high dimensions. A possible extension of   $\varphi$
 would have to take into account the minimal complexity 
of the branch locus,  (e.g. its Betti number) over all 
branched coverings. For a given  $N$ this complexity  
must be bounded from above, as it does happen 
in the case when  $N$ is a sphere by Alexander's theorem. 
However, it seems that such invariants are not easily computable.
\end{remark}

 \bibliographystyle{plain}

\end{document}